\documentclass{article}
\usepackage{amssymb,amsmath,ifsym}

\usepackage[mathscr]{euscript}
\usepackage{graphics}
\input{epsf}

\newfont{\weird}{cmff10}

\newcounter{sec}

\newcounter{punct}[sec]

\def\punct{\refstepcounter{punct}{\arabic{sec}.\arabic{punct}.  }}

\newtheorem{theorem}{Theorem}[sec]
\newtheorem{proposition}[theorem]{Proposition}

\newtheorem{lemma}[theorem]{Lemma}

\newtheorem{corollary}[theorem]{Corollary}

\newtheorem{question}[theorem]{Question}

\def\COUNTERS{\addtocounter{sec}{1}
              \setcounter{punct}{0}
          \setcounter{equation}{0}
          \setcounter{theorem}{0}
          }
          
          \def\sm{\smallskip}

\begin{document}
	

\newcommand{\rk}{\mathop {\mathrm {rk}}\nolimits}
\newcommand{\Aut}{\mathop {\mathrm {Aut}}\nolimits}
\newcommand{\Out}{\mathop {\mathrm {Out}}\nolimits}
\renewcommand{\Re}{\mathop {\mathrm {Re}}\nolimits}

\def\ov{\overline}
\def\un{\underline}

\def\wt{\widetilde}
\def\wh{\widehat}

\renewcommand{\rk}{\mathop {\mathrm {rk}}\nolimits}
\renewcommand{\Aut}{\mathop {\mathrm {Aut}}\nolimits}
\renewcommand{\Re}{\mathop {\mathrm {Re}}\nolimits}
\renewcommand{\Im}{\mathop {\mathrm {Im}}\nolimits}

\def\bfa{\mathbf a}
\def\bfb{\mathbf b}
\def\bfc{\mathbf c}
\def\bfd{\mathbf d}
\def\bfe{\mathbf e}
\def\bff{\mathbf f}
\def\bfg{\mathbf g}
\def\bfh{\mathbf h}
\def\bfi{\mathbf i}
\def\bfj{\mathbf j}
\def\bfk{\mathbf k}
\def\bfl{\mathbf l}
\def\bfm{\mathbf m}
\def\bfn{\mathbf n}
\def\bfo{\mathbf o}
\def\bfp{\mathbf p}
\def\bfq{\mathbf q}
\def\bfr{\mathbf r}
\def\bfs{\mathbf s}
\def\bft{\mathbf t}
\def\bfu{\mathbf u}
\def\bfv{\mathbf v}
\def\bfw{\mathbf w}
\def\bfx{\mathbf x}
\def\bfy{\mathbf y}
\def\bfz{\mathbf z}

\def\bfA{\mathbf A}
\def\bfB{\mathbf B}
\def\bfC{\mathbf C}
\def\bfD{\mathbf D}
\def\bfE{\mathbf E}
\def\bfF{\mathbf F}
\def\bfG{\mathbf G}
\def\bfH{\mathbf H}
\def\bfI{\mathbf I}
\def\bfJ{\mathbf J}
\def\bfK{\mathbf K}
\def\bfL{\mathbf L}
\def\bfM{\mathbf M}
\def\bfN{\mathbf N}
\def\bfO{\mathbf O}
\def\bfP{\mathbf P}
\def\bfQ{\mathbf Q}
\def\bfR{\mathbf R}
\def\bfS{\mathbf S}
\def\bfT{\mathbf T}
\def\bfU{\mathbf U}
\def\bfV{\mathbf V}
\def\bfW{\mathbf W}
\def\bfX{\mathbf X}
\def\bfY{\mathbf Y}
\def\bfZ{\mathbf Z}

\def\frD{\mathfrak D}
\def\frL{\mathfrak L}

\def\bfw{\mathbf w}

\def\R {{\mathbb R }}
 \def\C {{\mathbb C }}
  \def\Z{{\mathbb Z}}
  \def\H{{\mathbb H}}
\def\K{{\mathbb K}}
\def\N{{\mathbb N}}
\def\Q{{\mathbb Q}}
\def\A{{\mathbb A}}
\def\V{{\mathbb V}}
\def\O{{\mathbb O}}

\def\T{\mathbb T}
\def\P{\mathbb P}

\def\G{\mathbb G}

\def\cB{\EuScript B}
\def\cD{\EuScript D}
\def\cH{\mathscr H}
\def\cL{\mathscr L}
\def\cK{\EuScript K}
\def\cM{\EuScript M}
\def\cN{\EuScript N}
\def\cO{\EuScript O}
\def\cF{\EuScript F}
\def\cR{\EuScript R}
\def\cX{\EuScript X}

\def\bbA{\mathbb A}
\def\bbB{\mathbb B}
\def\bbD{\mathbb D}
\def\bbE{\mathbb E}
\def\bbF{\mathbb F}
\def\bbG{\mathbb G}
\def\bbI{\mathbb I}
\def\bbJ{\mathbb J}
\def\bbL{\mathbb L}
\def\bbM{\mathbb M}
\def\bbN{\mathbb N}
\def\bbO{\mathbb O}
\def\bbP{\mathbb P}
\def\bbQ{\mathbb Q}
\def\bbS{\mathbb S}
\def\bbT{\mathbb T}
\def\bbU{\mathbb U}
\def\bbV{\mathbb V}
\def\bbW{\mathbb W}
\def\bbX{\mathbb X}
\def\bbY{\mathbb Y}

\def\kappa{\varkappa}
\def\epsilon{\varepsilon}
\def\phi{\varphi}
\def\le{\leqslant}
\def\ge{\geqslant}

\def\B{\mathrm B}
\def\Ams{\mathrm{Ams}}

\def\la{\langle}
\def\ra{\rangle}

\def\F{{}_2F_1}
\def\FF{{}_2F_1^\C}

\newcommand{\dd}[1]{\,d\,{\overline{\overline{#1}}} }

\def\lambdA{{\boldsymbol{\lambda}}}
\def\alphA{{\boldsymbol{\alpha}}}
\def\betA{{\boldsymbol{\beta}}}
\def\mU{{\boldsymbol{\mu}}}
\def\sigmA{{\boldsymbol{\sigma}}}
\def\taU{{\boldsymbol{\tau}}}

\def\1{\boldsymbol{1}}
\def\2{\boldsymbol{2}}

\def\9{\textifsym{9}}
\def\8{\textifsym{8}}
\def\0{\textifsym{0}}

\begin{center}
	{\bf\large
	  Generalized eigenfunctions of
	  measure preserving transformations and  Riesz products}


\bigskip

\sc\large Yury A. Neretin%
\footnote{Supported by the grants FWF, P28421, P31591, PAT5335224,
and by MSHE RF GZ project.}
\end{center}

{\small We consider rank one measure preserving transformations $g$ and the corresponding Koopman
unitary operators $U(g)$.
It is known that  a generic (in the sense of Baire category) measure preserving  transformation has rank one,
  spectral
measure of $U(g)$ is purely singular and is given by a Riesz product.
For such transformations we  write explicitly 
spectral decompositions of operators $U(g)$ in generalized eigenfunctions.
}

\section{Three  questions\\ about measure preserving transformations}

\COUNTERS

{\bf \punct Spectral decompositions and integral transforms.}
There is a collection of second order self-adjoint  ordinary differential operators, for which decompositions in eigenfuctions
are explicit. There are  also such collections of difference operators on lattices, $q$-difference operators, difference operators in imaginary
direction (on the last type, see \cite{Ner-imaginary}).
Such operators with pure discrete spectra  produce the whole
Askey--Wilson hierarchy of hypergeometric orthogonal polynomials, i.e.,
Hermite, Laguerre--Sonine, Charlier, Jacobi (including Chebyshev and Gegenbauer),
Meixner, Meixner--Pollaczek,
Krawtchouk,  Hahn, dual Hahn, continuous Hahn,
 continuous dual Hahn, Racah, 
and
Wilson polynomials, see  \cite{Koe}.
Operators with continuous  spectra produce a collection of integral transforms, as the  Fourier transform, the Hankel transform,
the Mehler--Fock transform, the Kontorovich--Lebedev transform,
the Whittaker transform, 
the Olevsky--Jacobi transform, and many transforms, which do not have
well-recognized names.

 Such spectral decompositions appeared initially in classical theory
of boundary value problems and in elementary noncommutative harmonic analysis. Now they additionally   have numerous
applications in theory of special functions, see, for example,
\cite{Ner-wilson}, \cite{Koel}.

There also  finite systems of hypergeometric orthogonal polynomials
of Romanovsky type, see a collection in \cite{Koe} (it is not complete,
see for example \cite{Ner-wilson}). In fact, each Romanovsky system
is a discrete part of spectrum of an explicitly solvable Sturm--Liouville
problem with a mixed continuous-discrete spectrum.
There are also many solvable spectral problems for families of commuting operators
and  noncommutative spectral problems related to Plancherel formulas for
unitary representations.

In any case,  explicit spectral decompositions imply numerous
 consequences in  special functions, integral transforms, 
 noncommutative harmonic analysis, and mathematical physics.
 The topic of our paper are spectral decompositions of Koopman
 operators corresponding to measure preserving transformation. 
 
\sm

{\bf\punct Ergodic transformations and singular spectral measures.}
Consider the group $\Ams(M,\mu)$ of all measure preserving transformations 
of a space $(M,\mu)$ with finite or $\sigma$-finite non-atomic measure.
It is  a topological group equipped with a topology of Polish space, in particular
Baire category is well-defined on this group.
Ergodic transformations form a comeagre set.

Consider $g\in \Ams(M,\mu)$ and 
the corresponding Koopman unitary operator $U(g)$ in $L^2(M)$,
$$
U(g)\,f(m)=f(g(m))
$$
 (if $M$ is a space with a finite measure,
 then it is reasonable to consider
the action in the space of functions with zero mean).
It is known that spectral
measures of generic (in the sense of Baire category) measure preserving transformations are quite unusual from the point of view of the classical analysis.
For a dense $G_\delta$-set in $\Ams(M,\mu)$ a spectrum of $U(g)$ is multiplicity free,
spectral measures $\nu$ are purely singular (see \cite{Stepin}), and all convolution 
powers of $\nu$, i.e., $\nu$, $\nu*\nu$, $\nu*\nu*\nu$,  $\dots$ are pairwise
singular (see Katok, Stepin \cite{KS}, see also \cite{Nad2}).

\sm

 Notice that a similar statement holds
for generic unitary operators (see \cite{CN1}, \cite{Nad2}). Singular spectra are usual
for Schr\"odinger operators, see Chulaevsky \cite{Chu}, Simon \cite{Sim},
they appear for natural difference operators \cite{SiSp},
\cite{Jit}. 

\sm

{\bf \punct The questions.}

\begin{question}
Is it possible (in some interesting cases) to write explicitly spectral decompositions
for Koopman operators $U(g)$ with non-discrete spectra?
\end{question}

\begin{question}
Is it possible (in some cases) to write explicitly spectral decompositions
for Koopman operators with singular spectral measures?
\end{question}

\begin{question}
Is it possible to get nontrivial special function type identities from
inversion formulas and Plancherel formulas for measure preserving transformations?
\end{question}

The first question was formulated by S.~V.~Fomin \cite{Fom}, 1955. It seems
that his note did not entail any further continuations.
In the next section, we present definitely affirmative answers
 to Questions 1.1 and 1.2.
 
\sm

{\bf \punct Further structure of the paper.} 
In Subsect. \ref{ss:smooth}--\ref{ss:torus}, 
 we 
formulate several remarks (arising to Fomin) for a clarification of the structure
of the problem.

In Sections 2-3, we write explicit decompositions in eigenfunctions 
for operators corresponding to rank one transformations.

In Section 4, we play with spectral decomposition and construct an
orthogonal basis in $L^2$ with respect to a singular measure;
it consists of rational functions on the circle.

\sm

{\bf\punct Smooth dynamical systems.%
\label{ss:smooth}} Consider a smooth compact
$n$-di\-men\-sio\-nal
manifold $M$ equipped with a $C^\infty$-smooth probabilistic measure. Let $g$
be a $C^\infty$-smooth measure preserving diffeomorphism $g$.
Let $U(g)$ be the corresponding unitary
 operator in $L^2(M)$. Then we can apply the formalism of rigged spaces and generalized eigenfunctions 
invented by L.~G\aa rding and I.~M.~Gelfand--A.~G.~Kostychenko
(see different variants of its exposition in Gelfand, Vilenkin \cite{GV}, Berezin, Shubin \cite{BeS}, Berezansky, Us, Sheftel
\cite{BSU}). We consider the triple
\begin{equation}
C^\infty(M)\subset L^2(M)\subset \cD'(M),
\label{eq:ridging1}
\end{equation}
where $\cD'(M)$ denotes the space of distributions on $M$. The space $C^\infty(M)$ is nuclear and invariant with respect to $U(g)$. Therefore
(see, e.g., \cite{GV})
there is a system of generalized eigenfunctions 
$\Xi_\sigma(m)\in \cD'(M)$ enumerated by points $\sigma$ ranging
in some Lebesgue measure space $\Sigma$ (and measurably depending on
$\sigma$), such that 

\sm 

--- $U(g)\,\Xi_\sigma =\lambda(\sigma)\, \Xi_\sigma$,
where $\lambda(\sigma)\in \C$ and $|\lambda(\sigma)|=1$;
 
\sm

--- the operator 
$$
J f(\sigma)=\int_M f(m)\cdot \Xi_\sigma(m)
$$
 sending $C^\infty(M)$ to the space of measurable functions
 on $\Sigma$  admits an extension to a unitary operator
$L^2(M)\to L^2(\Sigma)$; here the symbol $\int$ denotes the pairing
$C^\infty(M)\times \cD'(M)\to\C$.

\sm

Moreover, we can replace the rigging \eqref{eq:ridging1}
by
\begin{equation}
H^{n/2+\epsilon}(M)\subset L^2(M)\subset H^{-n/2-\epsilon}(M),
\end{equation} 
where $H^s(M)$ is a Sobolev space in the standard notation, see \cite{BSU},
and $\epsilon>0$ is arbitrary
(since the tautological embeddings $H^{n/2+\epsilon}(M)\to L^2(M)$,
$L^2(M)\to H^{-n/2-\epsilon}(M)$ are Hilbert--Schmidt operators,
see \cite{BeS}, Suppl.1.2.3).

 \sm
 
 {\bf \punct An  experiment. Automorphisms  of torus.%
 \label{ss:torus}}
 Consider the torus $\T^2=\R^2/2\pi\,\Z^2$ with coordinates $(\phi,\psi)$
 defined modulo $2\pi$.
 Consider the transformation of the torus defined by the matrix
 $A=\begin{pmatrix}2&1\\1&1\end{pmatrix}$,
 $$
 \begin{pmatrix}\phi&\psi\end{pmatrix}\mapsto
   \begin{pmatrix}\phi&\psi\end{pmatrix} A
.$$
The transformations
 $A$ and $A^{-1}=\begin{pmatrix}1&-1\\-1&2 \end{pmatrix}$
 of $\R^2$  preserve the lattice $2\pi k+2\pi l$, where $l$, $k$ range in $\Z$. Therefore they are well defined on the torus; since $\det A=1$, they are
 area preserving. Denote by $H$ the group consisting
of transformations $A^k$, where $k\in\Z$. Notice that
$$
A^k=\begin{pmatrix}h_{2k+1}&h_{2k}\\h_{2k}&h_{2k-1}\end{pmatrix},
$$
where 
$$h_k:=
\frac1{\sqrt 5}
\Bigl[
\Bigl(\,\frac{1+\sqrt 5}{2}\,\Bigr)^k-
\Bigl(\,\frac{1-\sqrt 5}{2}\,\Bigr)^k\Bigr]
$$
are the Fibonacci numbers.

 We have the canonical identification 
 $$
 L^2(\T^2)\simeq \ell(\Z^2)
 $$ 
deined by the Fourier expansion. So, we come to the transformation
of the lattice $\Z^2$ defined by 
$\begin{pmatrix}m\\n\end{pmatrix}\mapsto
 A \begin{pmatrix}m\\n\end{pmatrix}$ and the corresponding unitary 
 operator in 
$\ell^2(\Z^2)$. Clearly, the group $H$ acts on  $\Z^2\setminus(0,0)$
 freely; this action preserves  'hyperbolas'
\begin{equation}
m^2+2mn-n^2=\rm{const}.
\label{eq:hyperbolas}
\end{equation}
 So $\ell^2(\Z^2)$ splits into a direct sum 
 $\oplus_j\ell^2(\cO_j)$, where $\cO_j$ are orbits of $H$ on 
 $\Z^2$ and subspaces $\ell_2(\cO_j)$ are $H$-invariant. 
 We have the trivial orbit $\cO_0$ consisting of one point 
  $(0,0)$. On remaining  orbits $\cO_j$ the group $H$ act freely, 
  each $\ell^2(\cO_j)$ 
   is $\ell^2(\Z)$, and the restriction
  of $U(A)$ to $\ell^2(\cO_j)$ is the usual  shift 
  operator. So, the problem of spectral decomposition
  reduces to a trivial question. For each $j$ we choose an element
  $\begin{pmatrix}m_j\\n_j\end{pmatrix}\in\cO_j$. For each pair
  $(j, \lambda)$, where $|\lambda|=1$, we write the distribution
  $$
  \Xi_{j,\lambda}(\phi,\psi):=
  \sum_{k=-\infty}^\infty \lambda^{-k} 
  \exp\biggl\{\begin{pmatrix} \phi&\psi\end{pmatrix}A^{k}
  \begin{pmatrix}m_j\\n_j\end{pmatrix} \biggr\}
  $$ 
 (the series converges as a series of distributions on the torus). 
 Also, $\xi(x,y)=1$ is a fixed vector.
Clearly, we get a complete collection of generalized eigenfunctions.
Nontrivial eigenfunctions are enumerated by points of $\Z\times \R$;
this space is
equipped with the product of
 the counting measure on $\Z$ and the Lebesgue measure $d\lambda$
on $\R$. 

Notice that the  spectrum has infinite multiplicity,
in classical theories  spectra of multiplicity $\ge 1$
 are serious  
obstacle for further applications.   

 \sm
 
  The matrix $A$ can be replaced by an arbitrary
$2\times 2$ matrix $B$ with integer coefficients such that $\det B=\pm1$
and eigenvalues of $B$ are noncoinciding reals.
Moreover, we can replace $\T^2$ by an arbitrary separable Abelian locally compact group $K$ and $A$ by an arbitrary ergodic automorphism $\theta$ of $K$.
The dual group $\wh K$ is discrete countable.
Applying the Fourier transform, we identify $L^2(K)$ with
$\ell^2(\wh K)$, the dual automorphism $\wh \theta$ is 
a permutation of $\wh K$.

\sm

However, in our trivial example
we observe unusual distributions $\Xi_{j,\lambda}$.
I can not say much about them.
Clearly, they are eigendistributions
of the operator
$$
\Delta:=\frac{\partial^2}{\partial \phi^2}+
\frac{\partial^2}{\partial \phi\, \partial \psi}-\frac{\partial^2}{\partial \psi^2}
$$
on the torus $\T^2$.
 It can be shown
  that  wave fronts of all such distributions consist of two
 lines $\eta=\pm (2-\sqrt 5)\,\xi$, where $(\xi,\eta)$
 are the coordinates dual to $(\phi,\psi)$;
   these lines are the asymptotes
  of hyperbolas of \eqref{eq:hyperbolas}. 
  
  Notice that the joint spectrum of commuting operators $U(A)$
  and $H$ has finite multiplicities.

\sm

\section{The statement}

\COUNTERS

{\bf \punct Purpose of the paper.}
A generic measure preserving transformation has rank one.
  Any  rank one transformation
can be represented in a transparent form
using cutting and stacking model,
 which  was introduced by Ornstein \cite{Orn},
 see for example \cite{Fri}, \cite{Nad2}
(for an equivalence of the abstract and constructive definitions, see \cite{Bax}).
However, a verification of conjugacy/non-conjugacy
of transformations having such 'normal forms' is a non-obvious question, see an interesting
discussion of a conjugacy problem in \cite{FRW}. 

Inspite of properties, which could be regarded as bad,
spectral measures for rank one transformations 
admit wonderful explicit expressions  in terms of
(generalized%
\footnote{For original F.~Riesz products (1918)
see \cite{Rie}, \cite{Zygmund}.}) Riesz products
$$
\prod_{j=1}^\infty \frac 1{m_j} P_j(z)\ov{P_j(z)},
$$
where $P_j(z)$ are 'sparse' polynomials of the form
$$
P_j(z)=1+\sum_{k=1}^{m_j-1} z^{r_{jk}}.
$$
Under certain conditions of 'sparseness',
such products  converge in the sense of weak convergence of measures.
Apparently, the first example in this spirit was obtained by Ledrappier
\cite{Led}, 1970. Queff\'elec showed that Riesz products 
can arise as spectral measures for substitutional
dynamical systems \cite{Que}.
Bourgain \cite{Bou}, 1993, discovered the general formulas for  spectral measures
for rank one transformations and also showed
that for a generic (in a certain sense) rank one transformation 
spectral measures 
are determined by Riesz products.
See also%
\footnote{Consider the direct sum $H$ and   the direct product $X$ of a countable number of $\Z_2$.
 Ismagilov in \cite{Ism1}-\cite{Ism2}, 1986, considered actions of $H$ in $L^2(X)$,
 'induced from a virtual subgroup' in the Mackey sense and 
obtained Riesz products as spectral measures for representations of this kind.}, \cite{CN}, \cite{KR}, \cite{Nad2}.

The purpose of this paper is to obtain explicit spectral decompositions
of operators $U(g)$ for rank one transformations $g$, see Theorem \ref{th:1}.

In Subsect. \ref{ss:ring}--\ref{ss:rank-one}, we present a convenient for our purposes description 
of rank one transformations.
The standard model with cutting and stacking intervals  is explained later in Subsect. \ref{ss:cutting}.
In that construction both  measure space and its transformation
are defined in inductive way. 
 Our description in Subsect. \ref{ss:ring}--\ref{ss:rank-one} is  based on 'coordinates' on a measure
space, this allow us to write formulas.

\sm

{\bf \punct Rings $\O[\bfm]$.%
\label{ss:ring}}
For natural $k$, $l$ we have a canonical homomorphism of residue rings
$$
\Z/\,l\Z\leftarrow  \Z/\,k l\Z.
$$
 Fix a sequence $\bfm=m_1$, $m_2$, \dots\vphantom{.} of integers,
$m_j\ge 2$.
 Consider the inverse limit
 $$
 \O=\O[\bfm]=\O[m_1,m_2,\dots]=\lim\limits_{\longleftarrow}
  \Z/m_1 m_2 m_2\dots m_k\Z
 $$
 of residue rings
\begin{equation}
\Z/m_1\Z\leftarrow \Z/m_1 m_2\,\Z\leftarrow  \Z/m_1 m_2 m_3\,\Z\leftarrow \dots
\label{eq:inverse}
\end{equation}
We write elements of $\O[\bfm]$ as 'numbers' infinite to the left
\begin{equation}
\bfs:=\dots s_3 s_2 s_1;\qquad \text{where $s_j$ ranges in the set $\{0,1, \dots, m_{j}-1\}$}.
\label{eq:bfs}
\end{equation}
For $\bfs \in\O[\bfm$] we assign element 
$\ov{\ov{s_j\dots s_1}}$ of $\Z/m_1\dots m_j\Z$
by
\begin{equation}
\ov{\ov{s_j\dots s_1}}:=s_1+ s_2 m_1+s_3 m_1m_2m_3+\dots +s_j m_1\dots m_{j-1}.
\label{eq:two-lines}
\end{equation}
and this determines the canonical surjective homomorphism
 $$\O[\bfm]\to \Z/m_1\dots m_j\Z.$$

We use the notation
$$
\9_{\un j}:=m_j-1,
\qquad
\8_{\un j}:=m_j-2.
$$ 
Also, we write 
$$
0_{\un j},\qquad 1_{\un j}
$$
if we want to show a position of $0$ or $1$ in the 'digital' notation (\ref{eq:bfs}).

Denote by $\zeta(\bfs)$ number of the first digit different from $\9$  at the end of $\bfs$, 
by $\iota(\bfs)$ this digit. Namely, for
\begin{equation}
\bfs:=\dots k\, \9_{\un j} \9_{\un{j-1}}\dots \9_{\un 1};\qquad k\ne \9_{\un{j+1}},
\label{eq:9}
\end{equation}
 we set
\begin{equation}
\zeta (\bfs):=j+1,\qquad \iota(\bfs):=k.
\label{eq:zeta}
\end{equation}
For the number $ \dots \9_{\un 3}\9_{\un 2}\9_{\un 1};$ these functions are not defined.

By the definition, $\O[\bfm]$ is a  ring, addition and product can be easily described
in terms of columnar addition and multiplication.
Below we need only  the operation
of addition of 1, it  can be defined in terms of columnar addition as follows. Consider
$$
+
\begin{array}{r}
\dots s_3 s_2 s_1;
\\
1;
\end{array}
$$
If $s_1\ne\9_{\un 1}$, then we simply replace $s_1\mapsto s_1+1$.
If
$\bfs$ has the form (\ref{eq:9})
then 
$$
\bfs+1= \dots (k+1) 0_{\un j} 0_{\un{j-1}}\dots 0_{\un 1}.
$$

\sm

The set $\O[\bfm]$ has a natural structure of a compact metric space.
Let $\bfs$, $\bft\in \O[\bfm]$ have the form
\begin{align*}
\bfs:=\dots k \,\sigma_j\dots\sigma_1;
\qquad
\bft:=\dots l\, \sigma_j\dots\sigma_1;
\end{align*}
and 
$l\ne k$. Then
$$
d(\bfs,\bft)= \frac 1{2^j} 
.
$$
This space is ultrametric. In particular, for any pair $B_1$, $B_2$
 of closed balls  we have $B_1\cap B_2=\varnothing$, or $B_1\supset B_2$,
 or $B_2\supset B_1$. We say that a {\it ball of level $k$} is a ball of  
 radius $2^{-k}$.
  Such balls $B[\sigmA]$ are enumerated by
 finite 'numbers'
 $$
  {\sigmA}:=\sigma_k \sigma_{k-1} \dots \sigma_1;
 $$
 the  corresponding ball consists of all sequences
 $$
 \dots s_{k+2} s_{k+1}\sigma_k \sigma_{k-1} \dots \sigma_1;
 $$
where $s_{k+1}$, $s_{k+2}$, $\dots$  are arbitrary, $s_i\in\bigl\{0,1, \dots, \9_{\un i}\bigr\}$.

We define a canonical probabilistic {\it measure} $\mu$ on $\O[\bfm]$ assuming that  measure of any 
ball of level $k$ is $(m_1\dots m_k)^{-1}$.

\sm

{\sc Remark.} If all $m_j$ are equal to a prime $p$, then 
$\bbO[\bfm]$ is the ring $\O_p$ of $p$-adic integers. 
In general case, decompose each $m_j$ into prime factors,
$m_j=\prod_p p^{\alpha_j(p)}$.
Then the chain \eqref{eq:inverse} can be represented in the form
$$
\dots \longleftarrow
 \bigoplus_p\, \Z\Bigl/\Bigl(\prod_{i\le j} p^{\alpha_i(p)} \cdot \Z\Bigr) \longleftarrow\dots
$$
Therefore,
$$
\lim\limits_{\longleftarrow} \Z/m_1\dots m_j \Z=
\bigoplus_p \lim\limits_{\longleftarrow} \Z\Bigl/
\Bigl( \prod_{i\le j}p^{\alpha_i(p)}\cdot \Z\Bigr).
$$
So,  $\O(m_1,m_2,\dots)$ is a product $\prod_p \bfO_p$
over all primes $p$; a factor $\bfO_p$ is the ring   $\bbO_p$
of $p$-adic integers
if $\sum_j \alpha_j(p)=\infty$, and a residue ring
$\Z/p^N \Z$ if $\sum_j \alpha_j(p)=N$.
For the $p$-adics integers $\O_p$ and integer adeles $\prod_p\O_p$,
 there is 
a commonly recognized Bruhat's definition of distributions, and this definition makes sense for arbitrary rings $\bbO[\bfm]$.
 Namely, we have the space
$\cF(\bbO)$ of locally constant functions, the space of distributions 
$\cF'(\bbO)$ is dual to  $\cF(\bbO)$.
Thus,
$$\cF(\bbO)\subset L^2(\bbO)\subset \cF'(\bbO).  $$
We emphasize that the construction below depends
on the sequence $m_j$ (and not only on the sums $\sum_j\alpha_j(p)$).
\hfill $\boxtimes$

\sm

{\bf \punct Rank one transformations.%
\label{ss:rank-one}}
Parameters of  rank one measure preserving transformations
are  following:

\sm

$\bullet$ We fix integers $m_j\ge 2$ as above.

\sm

$\bullet$
For each
$j=1$,  2, 3, \dots, we fix nonnegative integers
$a_j^0$, $a_j^1$, \dots, $a_j^{m_j-2}$. 

\sm

For such data we define a locally compact totally disconnected topological space%
\footnote{This space is metrizable (as a union of a countable number of disjoint clopen compact sets),
but we do not define a metric.}
$\V=\V[\bfm,\bfa]$ and a transformation
$Q=Q[\bfm,\bfa]$ of $\V[\bfm,\bfa]$ as follows. 

\sm

{\sc The space $\V=\V[\bfm,\bfa]$.}
 We consider 'numbers' in $\bbO[\bfm]$ having one additional  digit after 
the semicolon,
$$
\bfs:=\dots s_3 s_2 s_1; s_0,
$$
where $\dots s_3 s_2 s_1\in \O$ are the same as above,
$$
\dots s_3 s_2 s_1;\ne \dots \9_{\un 3} \9_{\un 2} \9_{\un 1};
$$
(this point is excluded)
and  
\begin{equation}
 s_0\in \bigl\{0,1, \dots, a_{\zeta(\bfs)}^{\iota(\bfs)}\bigr\}.
 \label{eq:s0}
\end{equation}

Let us repeat this less formally. Take $\dots s_3 s_2 s_1\in \O[\bfm]$.
We look to the end of this 'number',
$$
\dots
\iota\, \9_{\un{\zeta}-1}\dots \9_{\un 1};
$$
where
$$
\zeta=\zeta(\dots s_3 s_2 s_1;), \qquad \iota:=\iota(\dots s_3 s_2 s_1;)
$$
Then a digit after the semicolon ranges in the set  (\ref{eq:s0}).

\sm

Denote by $\V:=\V[\bfm,\bfa]$ the set of all such sequences.
We denote  $a_{\zeta}^{\iota}$ after the semicolon by $\9^{\un {\zeta,\iota}}_{\un 0}$.

\sm

{\sc Balls.}
Next, we define {\it balls $B_k[\sigmA]$ of level $k$}.
We fix a finite sequence 
 $$
 \sigmA:= \sigma_k \sigma_{k-1} \dots \sigma_1;\sigma_0
,$$
where $\sigma_j\in \{0,1,\dots,m_j-1\}$ are as above, we only require
\begin{equation}
\sigmA^\circ:=\sigma_k \sigma_{k-1} \dots \sigma_1;\ne \9_{\un k}\9_{\un{k-1}}\dots \9_{\un 1};
\label{eq:excl}
\end{equation}
and
 $$\sigma_0\in\bigl\{0,\dots,
  a_{\zeta(\sigmA^\circ)}^{\iota(\sigmA^\circ)}\bigr\}.$$
  Now we
are ready to define the ball  $B_k[\sigmA]$. It consists of all 'numbers'
$$
\dots s_{k+2} s_{k+1}\sigma_k \sigma_{k-1} \dots \sigma_1;\sigma_0
.
$$

\sm

 There are $m_1-1$ balls $B_1[j]$ of level $1$. Such ball consists
 of elements of the form
 $$
 \dots s_4 s_{3} \dots s_2 j_1;0,
 $$
 where $j=0$, 1, \dots, $m_j-2$ is fixed (recall that we assume that $a_1^l=0$).

There are no balls of level 0. However, below we use the notation
$B_0[\mathbf{0}]:=\O$.

\sm

Each ball $I_k[\sigmA]$ admits a canonical division
into $m_{k+1}$ disjoint subballs of level $k+1$:
\begin{equation}
B_k[\sigma_k \sigma_{k-1} \dots \sigma_1;\sigma_0]=
\coprod_{\tau=0}^{m_{k+1}-1}
B_{k+1}[\tau \sigma_k \sigma_{k-1} \dots \sigma_1;\sigma_0].
\label{eq:subballs}
\end{equation}


Denote by $h_k^\circ$ the number of balls of level $k$.
Denote
$$h_k=h_k^\circ+1,$$
 these numbers satisfy the recurrence
relation
\begin{equation}
h_{k}=m_{k} h_{k-1} +\sum_{i=0}^{m_k-2} a_{k}^{i}, \qquad h_0=1, \quad
h_1=m_1-1.
\label{eq:h_k}
\end{equation}

{\sc Remark.} Due to the exclusion (\ref{eq:excl}),
balls of level $k$ do not cover the whole space $\V[\bfm,\bfa]$.
\hfill $\boxtimes$

\sm

{\sc The measure on $\V[\bfm,\bfa]$.}
To define the measure $\mu$ on $\V[\bfm,\bfa]$, we
assign  measure $(m_1\dots m_k)^{-1}$ to each ball of level $k$.
The total measure of $\V[\bfm,\bfa]$ is 
$$
1+\sum_{k=1}^\infty \frac 1 {m_1 \dots m_k} \sum_{i=0}^{m_k-2} a_{k}^i.
$$
We admit both finite and infinite total measures.

\sm

{\sc The rank one transformation $Q$.}
We have a well-defined shift 
$$
Q:\bfs \mapsto \bfs+0;\!1
.
$$ 
Namely, if the last digit $\sigma_0$ of $\bfs\in \V[\bfm,\bfa]$ is not $\9_{\un 0}^{\un{\zeta,\iota}}$,
we simply replace $\sigma_0\mapsto \sigma_0+1$. If $\bfs$ has the form
$$
\bfs= \dots k \,\9_{\un j}\dots \9_{\un 1}; \9_{\un 0}^{\un{j+1,k}}, \qquad \text{where $k\ne \9_{\un {j+1}}$}
,$$
then 
$$
Q\bfs = \dots (k+1)\, 0_{\un j}\dots 0_{\un 1};0_{\un 0}^{\un{1,0}}.
$$
Thus we get a measure preserving map.

\sm

{\sc The inverse transformation.}
 We have also the transformation
$Q^{-1}:\sigmA\mapsto \sigmA-0;\!1$
defined for all $\sigmA\ne \dots 000;0$. Namely, if
the last digit $\sigma_0\ne 0$, then we simply replace 
$\sigma_0\mapsto\sigma_0-1$. If
$$
\sigmA=\dots k 0_{\un j}\dots 0_{\un 2}0_{\un 1};0_{\un 0}^{\un{j+1,k}},
\qquad\text{where $k\ne 0$,}
$$
then
$$
Q^{-1} \sigmA=\dots (k-1)\9_{\un j}\dots \9_{\un 2}\9_{\un 1};
\9_{\un 0}^{\un{j+1,k-1}}.
$$ 


\sm

{\bf \punct Test functions and distributions.}
For each ball $B_k[\sigmA]$ of level $k$
 we define its indicator function
 $$
 I_k[\sigmA](\bfs):=\begin{cases}
 1,&\qquad \text{ if $\bfs \in B[\sigmA]$;}
 \\
 0,&\qquad \text{ if $\bfs \notin B[\sigmA]$.} 
 \end{cases}
 $$
Clearly, we have
\begin{equation}
I_k[\sigma_k \sigma_{k-1} \dots \sigma_1;\sigma_0]=
\sum_{\tau=0}^{m_{k+1}-1}
I_{k+1}[\tau \sigma_k \sigma_{k-1} \dots \sigma_1;\sigma_0].
\label{eq:subdivision}
\end{equation} 
 
A {\it test function} is a finite linear combination
of indicator functions of balls. Without loss of generality,
we can assume that balls are disjoint. Otherwise, we pass
to a subdivision into smaller balls applying \eqref{eq:subdivision}.
We denote by $\cF(\V)$ the space of all test functions,
by $\cF'(\V)$ it dual space. Then we get a triple
\begin{equation}
\cF(\V)\subset L^2(\V)\subset \cF'(\V).
\label{eq:triple}
\end{equation}
We have the first embedding by definition.
Any element of $q\in L^2(\V)$ determines a linear functional
$$\kappa_q(h)=\int_\V q(\sigmA) h(\sigmA)\,d\sigmA$$
 on $\cF(\V)$.
So, we have the second embedding.

\sm

{\sc Remark.} The space $\cF(\V)$ is a nuclear space in the sense of 
the book Berezansky, Us, Shehtel, \cite{BSU}, Sect.14.2.2 (these
authors allow uncountably normed spaces). More precisely, denote by 
$\cF_k(\V)\subset \cF(\V)$ the subspace consisting of linear combination
of indicator functions of balls of level $k$. So,
$\cF_k(\V)\subset \cF_{k+1}(\V)$. Denote by $\cH_{k+1}(\V)$
the orthocomplement of $\cF_k(\V)$ in $\cF_{k+1}(\V)$.
Thus, any element $h$ of $\cF(\V)$ can be in a unique way represented 
as
$$
h=h_1+h_2+\dots,
$$ 
where $h_j\in \cH_j(\V)$ and $h_l=0$  for sufficiently large $l$.
For each sequence $C_j\ge 1$ we define the inner product   
$\la\cdot,\cdot\ra_C$ on  
$\cF(\V)$ by
$$
\la h,h'\ra_C:=\sum_{j=1}^\infty C_j\la h_j,h_j'\ra_{L^2(\V)}
$$  
Denote by $\ov{\cF}_C(\V)$ the completion of  $\cF(\V)$
by this inner product. Clearly, for any sequence $C_j$ there is
a sequence $D_j\ge C_j$ such that the identical embedding
$\ov{\cF}_D(\V)\to \ov{\cF}_C(\V)$ is a Hilbert--Schmidt operator. 
This means that ${\cF}(\V)$ is a nuclear space in the sense of 
\cite{BSU}. Now we can to apply the formalism of generalized eigenfunctions
to the triple \eqref{eq:triple}, formally we can refer to
 \cite{BSU}, Sect.15.1.3. 
\hfill $\boxtimes$ 
 
 \sm

A linear functional $\Xi\in \cF'(\V)$
is uniquely determined by its values 
$\xi_k[\sigmA]$ on indicator functions $I_k[\sigmA]$.
Clearly, such function must satisfy the following 
relations
$$
\xi_k[\sigma_k\dots\sigma_1;\sigma_0]
=\sum_{\tau=0}^{m_{k+1}-1} 
\xi_{k+1}[\sigma_k\dots\sigma_1;\sigma_0].
$$
Below we describe eigendistributions corresponding a given 
eigenvalue $z$ in terms of these numbers.
Equivalently, we present the formula 
\eqref{eq:phi} for  images of indicator functions in
the
spectral decomposition.

Formally, we do not use  abstract theorems about generalized eigenfunctions in our proof.

\sm

{\bf \punct The statement.} 
 Our purpose it to write an explicit spectral decomposition of
 the operator
 $$ U f(\bfs):=f(Q
 \bfs)$$
 in $L^2(\V)$.
 
 We need more notation. 
 For each ball 
 $B_k[\sigmA]=B_k[\sigma_k\dots \sigma_1;\sigma_0]$ of level $k$
 we define integers
 \begin{equation}
 \ov{\ov{\sigma_k \dots \sigma_{\beta}}}:=\sigma_{\beta}+\sigma_{\beta+1} m_{\beta+1}+
 \dots + \sigma_k m_{\beta+1}\dots m_k,\qquad \beta\ge 1,
 \label{eq:ov}
 \end{equation}
 and the function $\Upsilon_k$ on the set of balls of level $k$ by
 \begin{equation}
 \Upsilon_k(\sigmA):= \sum_{\alpha=1}^k
 \Bigl( \ov{\ov{\sigma_k \dots \sigma_{\alpha+1}}}\,
 \sum_{i} (a_{\alpha}^i+1) +\sum_{i<\sigma_{\alpha}}  (a_{\alpha}^i+1)
 \Bigr)+\sigma_0.
 \label{eq:Upsilon}
 \end{equation}
 
 Denote by $S^1$ the unit circle $|z|=1$  on the complex plane.
  Denote by  $d \dot z$ the standard probabilistic Lebesgue measure 
 on  $S^1$.
 For each $k$ we define a polynomial function $\Theta_k(z)$ on $S^1$ by
 \begin{equation}
 \Theta_k(z):= \sum_{p=0}^{m_k-1} z^{p h_{k-1}+\sum_{i<p} a_{k-1}^i }.
 \end{equation}
 
 Consider the sequence of measures
 \begin{equation}
 d\kappa_n:=
 \prod_{k=1}^n \frac 1{m_k} \Theta_k(z)\ov{\Theta_k(z)}\,
 \,d \dot z
  \label{eq:R1}
 \end{equation}
 and the weak limit (the {\it Riesz product})
 \begin{equation}
 d\kappa(z)=
 \prod_{k=1}^\infty \frac 1{m_k} \Theta_k(z)\ov{\Theta_k(z)}\,\, d \dot z:=\lim_{n\to\infty} d\kappa_n
 .
 \label{eq:R2}
 \end{equation}

Consider the map 
$$\cR:\cF(\V)\,\,\to\,\,\C^\infty(S^1)$$
  that sends each indicator function
$I_k[\sigmA]$ to the function $\Phi_k[\sigmA]$ on $S^1$
given by
\begin{equation}
\Phi_k[\sigmA]:= \frac{z^{\Upsilon_k(\sigmA)}}{\prod_{j=1}^k \Theta_k(z)}
.
\label{eq:phi}
\end{equation}

In particular, $\Psi_0[\mathbf{0}]=1$.

\begin{theorem}
	\label{th:1}
	If functions $\Theta_k(z)$ have no zeros on the circle, then the map
$\cR$ extends  to a unitary operator $L^2(\V)\to L^2(S^1,d\kappa)$
and
$$
\cR \,U f(z)=z\, \cR(z).
$$
\end{theorem}

\sm

{\bf \punct The correspondence with the cutting--stacking construction.%
\label{ss:cutting}} Recall  usual Ornstein's construction of rank one transformations.
Consider a sequence $m_1$, $m_2$, \dots, where $m_j\ge 2$. For each $j$
consider a collection $a_j^0$, \dots, $a_j^{m-1}$ of nonnegative integers.
Having this data, we construct a measure space $\cX$ and its measure
preserving transformation $\theta$ in the following inductive way. 

\sm

{\sc Preliminary comments.}
A measure space $\cX$ is a union
of an increasing chain
of spaces
$$
X_0\subset X_1\subset X_2\subset\dots.
$$
 Each space $X_k$ is equipped with   a measure
$\nu_k$ and  partition%
\footnote{We consider elements of a finite measurable partition
upto  perturbations by sets of zero measure.}
  $\xi_k$ into $h_k$ disjoint
pieces of measure $1/m_1\dots m_k$ (the number $h_k$
is defined in the process of the construction).
One element $A_k$ of the partition is fixed.
 The space $X_k$
is also equipped with
  a measure  preserving map
$\theta_k:X_k\setminus A_k\to X_k$;
this map sends elements of $\xi_k$ to elements of $\xi_k$.
Precisely, there is a set $B_k$ such that the partition $\xi_k$
consists of subsets $B_k$, $\theta_k B_k$, $\theta_k^2 B_k$,\dots,
$\theta_k^{h_k-1} B_k$, and the last set coincides with $A_k$.

 The embedding $X_k\to X_{k+1}$
is  measure preserving, and  $\theta_{k+1}\Bigr|_{X_k}=\theta_k$.
 Each element of
$\xi_k$ is a union of $m_{k+1}$ elements of the partition $\xi_{k+1}$.

 \sm

{\sc The base of the inductive construction.} For $k=0$ we assume
that $X_0$ is the segment $[0,1]$ and $A_0=X_0$.
Notice that $\nu_0(A_0)=\prod_{j=1}^0 m_j=1$. We assume $h_0=1$.

\sm

{\sc The inductive step.} After $k$-th step we have a disjoint union of $h_k$
segments of length $1/m_1\dots m_k$. We regard them  as pieces
of the segment 
$[0,h_k/m_1\dots m_k],$
this is our partition $\xi_k$.
 The set $A_k$ is the last
subsegment. The map $\theta_k$ is the shift 
$x\mapsto x+1/m_1\dots m_k$.
It sends pieces of the  partition $\xi_k$ as it is shown on Fig.~1.
The set $A_k$ is the last subsegment.

Next, we divide each subsegment 
$$
\Bigl[\frac n{m_1\dots m_k},\frac{n+1}{m_1\dots m_k}\Bigr]
$$
into $m_{j+1}$ equal pieces, see Fig.~2. These pieces are elements of the partition $\xi_{k+1}$ of $X_k$ that are contained in $\cX_{k}$.

The new space $X_{k+1}$ is $X_k$
with a union of $\sum_{i=0}^{m_{j+1}-1} a_{k+1}^i$  disjoint segments
of length $1/m_1\dots m_k m_{k+1}$.
It is convenient to draw these segments as $m_{j+1}$ columns over
$A_k$. Namely, we enumerate pieces of $A_k$ by $0$, $1$, \dots, $m_{j+1}$,
and draw $a_{k+1}^i$ segments over $i$-piece. These new
segments are elements of $\xi_{k+1}$. See Fig.~3.

It remains to extend the transformation $\theta_k$.  The map
$\theta_{k+1}$ sends each element of each column in the upper direction
by a shift (until achieving of the top). For $i\ne m_{k+1}-1$, we send
the top segment of $i$-th column by a shift to the $(i+1)$-th subsegment
of $[0,1/m_1\dots m_k]$.   The top of the last column is $A_{k+1}$.
See Fig.~4.

\begin{figure}
\epsfbox{stack.1}
\caption{A set $X_k$, its partition $\xi_k$, and the map $\theta_k$; here $h_k=5$.}
\end{figure}

\begin{figure}
\epsfbox{stack.2}
\caption{A subpartition of $X_k$; here $m_{k+1}=4$.}
\end{figure}

\begin{figure}
\epsfbox{stack.3}
\caption{The space $X_{k+1}$ and its partition $\xi_{k+1}$.
 Here
$a_{k+1}^0=2$, $a_{k+1}^1=1$, $a_{k+1}^2=4$, $a_{k+1}^3=2$.}
\end{figure}

\begin{figure}
\epsfbox{stack.4}
\caption{The map $\theta_{k+1}$ outside $X_k\setminus A_k$.}
\end{figure}

Now we have the space $X_{k+1}$ and its partition $\xi_{k+1}$ into
$$
h_{k+1}:=h_k\cdot m_{k+1}+\sum_{i=0}^{m_{k+1}-1} a_{k+1}^i
$$
segments of length $1/m_1\dots m_{k+1}$. Denote $J$ the most left
segment. Then the chain
$$
J,\,\theta_{k+1} J,\, \theta_{k+1}^2 J,\dots, \theta_{k+1}^{h_{k+1}} J
$$
contains all elements of the partition $\xi_{k+1}$, the last element is $A_{k+1}$.

This completes the inductive step.

\sm

{\sc The correspondence with the construction of Subsect.~\ref{ss:rank-one}.} First, let us assume that all 
$$a_k^{m_k-1}=0,$$
i.e., the last columns have high 0.

For each segment of level $l$ (i.e., element of $\xi_l$) we assign a ball of level $l$ in inductive way.

Assume that this is done for all balls of levels $\le k$.
Let an element of the partition $\xi_k$ different from $A_k$ corresponds
to a ball $B_k[\sigma_k\dots\sigma_1;\sigma_0]$.
Then its subsegments correspond to balls
$B_{k+1}[0\,\sigma_k\dots\sigma_1;\sigma_0]$, $B_{k+1}[1\,\sigma_k\dots\sigma_1;\sigma_0]$,
 $B_{k+1}[\9\sigma_k\dots\sigma_1;\sigma_0]$.
For an element of $i$-th column (where $i\ne m_{j+1}-1$) of floor $l$
we assign the ball
$B_{k+1}[i\, \9_k\dots\9_1;l]$.

It is easy to see that we get an isomorphism of the Boolean algebra generated by segments in $\cX$ and the Boolean algebra generated by balls
in $\V$.
Completing these algebras by the distances 
$$\rho_\cX(A,B):=\nu(A\bigtriangleup B) \quad\text{and}\quad
\rho_\V(S,T):=\mu(S\bigtriangleup T)
$$
(where $\bigtriangleup$ denotes the symmetric difference of sets)
we get an isomorphism of Boolean algebras of measurable sets. 
 This implies a point isomorphism of measure
spaces defined by a.s.
(see, for example  \cite{Bog}, Theorem 9.5.1).

\sm

Now, let $a_{k+1}^{m_{k+1}-1}\ne 0$. 
We must build column of this height  over 
$$
\Bigl[\frac{h_k-1}{m_1\dots m_k},\frac{h_k}{m_1\dots m_k}\Bigr].
$$
We do not do this, but on the next step
we change parameters  
$$
a_{k+2}^j\mapsto a_{k+2}^j +
a_{k+1}^{m_{k+1}-1},
$$ 
 etc.
 



\section{Proof of Theorem \ref{th:1}}

\COUNTERS

The proof given below is a straightforward verication
of the statement. Formally, it does not use any preliminary information
about rank one transformations, Riesz products, and spectral theory.

\sm

{\bf \punct The function $\Upsilon_k$.}
The map $Q$ sends any ball $B_k[\sigmA]$ of level $k$ to 
a ball of level $k$, namely to $B_k[\sigmA+0;\!1]$,  except the ball
\begin{equation}
B_k\bigl[\9_{\un k} \dots \9_{\un 2} \8_{\un 1}; \9_{\un 0}^{\un{1,m_1-2}}\bigr].
\label{eq:B9}
\end{equation}

\begin{lemma}
	\label{l:Upsilon}
	$$
B_k[\sigma_k\dots \sigma_1;\sigma_0]= Q^{\Upsilon(\sigma_k\dots \sigma_1;\sigma_0)}	
B_k[0_{\un k}\dots 0_{\un 1};0].
$$
\end{lemma}

{\sc Proof.}
We have a natural linear order on the set of balls of level $k$:
\begin{equation}
 B_k[\tau_k\dots \tau_1;\tau_0]
\prec
B_k[\sigma_k\dots \sigma_1;\sigma_0]
\label{eq:B-B-more}
\end{equation}
if for some $j$ we have
$\sigma_i=\tau_i$ for all $i>j$  and
$
\sigma_j>\tau_j
$.
If we have (\ref{eq:B-B-more}), then
$$
B_k[\sigma_k\dots \sigma_1;\sigma_0]= Q^L B_k[\tau_k\dots \tau_1;\tau_0]
$$
for some positive $L$. In particular,
$$
B_k[\sigma_k\dots \sigma_1;\sigma_0]=Q^N B_k[0_{\un k}\dots 0_{\un 1};0]
$$
for some positive  $N$. Clearly, $N$ is the cardinality of the set
of all balls $\prec B_k[\sigma_k\dots \sigma_1;\sigma_0]$.
Each ball of this type has the form
\begin{equation}
B[\tau_k\dots \tau_{\alpha+1}\, l\, \9_{\un{\alpha-1}}\dots \9_{\un 1},\tau_0],
\qquad \text{where $l\ne \9_{\un{m_\alpha}}$}
\label{eq:BBB}
\end{equation}
for some $\alpha$.

Thus, we fix $\alpha$ and evaluate number of balls (\ref{eq:BBB}),
which are $\prec B_k[\sigma_k\dots \sigma_1;\sigma_0]$.
This can happened in the following cases:

\sm

1) If
$$
\ov{\ov{\tau_k\dots \tau_{\alpha+1} }}
<\ov{\ov{\sigma_k\dots \sigma_{\alpha+1} }}.
$$
For each $\tau_k\dots \tau_{\alpha+1}$ we have 
$\sum_i (a_{m_\alpha}+1)$ balls
of this type.

\sm

2) If
$$
\ov{\ov{\tau_k\dots \tau_{\alpha+1} }}=
\ov{\ov{\sigma_k\dots \sigma_{\alpha+1} }}, \qquad l<\sigma_\alpha.
$$
Number of such balls  is $\sum_{i< \sigma_\alpha} (a_{\alpha}^i+1)$.

\sm

3) The exceptional case $\alpha=0$,
\begin{multline*}
\tau_k\dots \tau_{\beta}=\sigma_k\dots \sigma_{\beta}, \qquad
\sigma_\beta\ne \9_{\un\beta},\\
\tau_{\beta-1}\dots \tau_1=\sigma_{\beta-1}\dots \sigma_1=
\9_{\un{\beta-1}}\dots \9_{\un 1},\qquad \tau_0<\sigma_0.
\end{multline*}
Then we get $\sigma_0$ balls.

\sm

This gives us formula (\ref{eq:Upsilon}) for the exponent $N$.
\hfill $\square$

\sm

\begin{corollary}
	\label{cor:}
{\rm a)}
 Let $B_k[\sigmA]$ be a ball of level $k$ different from
 {\rm (\ref{eq:B9})}.
Then 
$$
\cR \,I_k[Q\sigmA]=z\cdot \cR\,I_k[\sigmA].
$$	

{\rm b)}  Let $B_k[\sigmA]$ be a ball of level $k$ different from
\begin{equation}
B\bigl[\dots 0_{\un k} 0_{\un j}\dots 0_{\un 2}0_{\un 1};0_{\un 0}^{\un{j+1,k}}\bigr].
\end{equation}
 Then 
 $$
 \cR I_k[Q^{-1}\sigmA]=z^{-1}\cdot\cR\, I_k[\sigmA].
 $$
\end{corollary}

We also get the following expression  for number  $h_k^\circ$ of balls of level $k$:
$$
h_k^\circ=\Upsilon_k\bigl(\9_{\un k}\dots \9_{\un 2} \8_{\un 1}; \9_{\un 0}^{\un{1,m_1-2}}\bigr).
$$

\begin{lemma}
The numbers $h_k$ satisfy the recurrence relation {\rm (\ref{eq:h_k})}.
\end{lemma}

{\sc Proof.}
Each ball of level $k$ is a union of $m_{k+1}$ subbals of level $(k+1)$.
This gives us $m_{k+1} (h_{k+1}-1)$ balls. Additionally,
we have balls of the type
$$
B_{k+1}\bigl[j\, \9_{\un k} \dots \9_{\un 1}; \tau \bigr].
$$
Number of such balls is 
$$\sum_{j=0}^{m_{k+1}-2} (a_{k+1}^j+1)=\sum_{j=0}^{m_{k+1}-2} a_{k+1}^j 
+m_{k+1}-1.$$
This gives our statement.
\hfill $\square$

\sm

{\bf \punct Self-consistency of the definition of the map $\cR$.}
Let us represent a ball  $B_{k}[\sigmA]$ of level $k$
as a union of $m_{k+1}$ disjoint balls of level $(k+1)$,
$$
B_{k}[\sigma_k\dots \sigma_1;\sigma_0]=
\coprod_{i=0}^{m_{k+1}-1}B_{k+1}[i\,\sigma_k\dots \sigma_1;\sigma_0].
$$
We must show that
$$
\cR\, I_{k}[\sigma_k\dots \sigma_1;\sigma_0]=
\sum_{i=0}^{m_{k+1}-1} \cR \, I_{k+1}[i\,\sigma_k\dots \sigma_1;\sigma_0].
$$
By Corollary \ref{cor:}, it is sufficient to consider the case
$$
\sigmA={0}_{\un k} \dots {0}_{\un 1}; {0},
$$
i.e., we must verify the following identity:
$$
\cR\, I_k[{0}_{\un k} \dots {0}_{\un 1};{0}]=
\sum_{j=0}^{m_{k+1}-1}\cR\, I_{k+1}[j_{k+1}{0}_{\un k} \dots {0}_{1};{0}].
$$
In the right-hand side we have
\begin{multline*}
 \sum_{j=0}^{m_{k+1}-1}\cR\,\, U^{jh_k+ \sum_{i<j} a_{k+1}^i} \,\,
 I_{k+1}[0_{\un{k+1}}{0}_{\un k} \dots {0}_{\un 1};{0}]
 =\\= \sum_{j=0}^{m_{k+1}-1} \frac{z^{jh_k+ \sum_{i<j} a_{k+1}^i}} {\prod_{p=1}^{k+1} \Theta_p(z)}
 =
 \frac {\Theta_{k+1}(z)}{\prod_{p=1}^{k+1} \Theta_p(z)}
 =
 \cR\, I_k[{0}_{\un k} \dots {0}_{\un 1};{0}].
 \qquad\square
 .
\end{multline*}

{\bf \punct The Riesz products.} 
Decompose
$$
P_k:=
\prod_{j=1}^k \frac 1{m_j} \Theta_j(z)\,\ov{\Theta_j(z)}=\sum_{\alpha=-N}^N
u_\alpha^{(k)} z^\alpha.
$$

\begin{lemma}
	\label{eq:Riesz}
{\rm a)} Coefficients $u_j^{(k)}$ satisfy
$$
0\le u_\alpha^{(1)}\le u_\alpha^{(2)}\le \dots \le 1.
$$

{\rm b)} $u_0^{(k)}=1$.

\sm

{\rm c)} In the product 
\begin{equation}
\prod_{j=l+1}^{k} \frac 1{m_j} \Theta_j(z)\,\ov{\Theta_j(z)}=\sum_{\alpha=-N}^N
u_\alpha^{(l:k)} z^\alpha
,
\label{eq:lk+1}
\end{equation}
we have 
$$
u_\alpha^{(l:k)}=0\qquad \text{fo all $\alpha$ satisfying $0<|\alpha|< h_l$.}
$$
\end{lemma}

{\sc Proof.}
We have
\begin{multline}
T_j(z):=
\frac 1{m_j} \Theta_j(z)\,\ov{\Theta_j(z)}
=\\=
1 + \frac 1{m_j}\sum_{0\le p<q<m_j}\Bigl( z^{(q-p) h_{j-1}+ \sum_{i:\, p\le i< q} a^i_{j-1} }
+\\+ z^{-\bigl((q-p) h_{j-1}+ \sum_{i:\, p\le i< q} a^i_{j-1}\bigr) }
\Bigr).
\label{eq:TT}
\end{multline}
The coefficients are nonnegative and this implies non-negativity and increasing 
of $u_\alpha^{(k)}$ (this is a part of the statement a).

The maximal degree of a monomial in (\ref{eq:TT})
 is 
$$
\nu_j:=(m_j-1)h_{j-1}  +\sum_{i=0}^{m_j-2} a_j^i.
$$ 
This value 
corresponds for $q=m_j-1$, $p=0$. 
We have 
\begin{equation}
\nu_j=h_j-h_{j-1}
\label{eq:nu-j}
\end{equation}
Therefore the maximal degree of a monomial in $P_k$ is
$$
\nu_1+\nu_2+\dots+\nu_k=(h_1-1)+(h_2-h_1)+\dots+(h_k-h_{k-1})=h_k-1,
$$
and minimal (negative) degree is $-h_k+1$.

On the other hand, minimal positive degree in 
$T_{k+1}$
is $h_{k}+a_{k+1}^0$,
it is achieved for $q=1$, $p=0$.
Therefore the term $z^0$ in $P_{k} T_{k+1}$ can arise only as a product of
$z^0\cdot z^0$. This implies statement b.

\sm

Thus $\int P_k(z)\,d\dot z=1$, the expression $P_k(z)$ is positive, therefore
other Fourier coefficients are $\le 1$. This completes a proof of the statement a.

\sm

Consider the product (\ref{eq:lk+1}), i.e.,
$
\prod_{j=l+1}^k T_j
$.
We expand each $T_j$ as (\ref{eq:TT}). Each term of the expansion has the form
\begin{equation}
z^{\xi_{j_1}} z^{\xi_{j_2}}\dots z^{\xi_{j_p}},
\end{equation}
where $z^{\xi_{j_s}}$ appears from the expansion of $T_{j_s}$,
\begin{equation}
l+1\le j_1<j_2<\dots <j_p\le k,
\label{eq:zxi}
\end{equation}
and $\xi_{j_s}\ne 0$ for all $s$ 
(i.e., we exclude  factors $z^0$ from (\ref{eq:zxi}),
for the term $\prod_{j=l+1}^k z^0$ expression  (\ref{eq:zxi})
is empty,
 empty product is 1). To be definite, assume that $\xi_{j_p}>0$.
Then $\xi_{j_p}\ge h_{j_p-1}$. On the over hand for $s<p$
we have $\xi_{j_s}\ge -\nu_{j_s}$. Therefore the total degree of a monomial can be estimated as
\begin{multline*}
\xi_{j_p}+\sum_{s<p} \xi_{j_s}\ge h_{j_p-1}-\sum_{s<p} 
\nu_{j_s}\ge h_{j_p-1}-\sum_{l+1<i<j_p} \nu_i=\\=
h_{j_p-1}-\sum_{l+1<i<j_p} (h_{i}-h_{i-1})=h_l.
\end{multline*}
This proves the last statement.
\hfill
$\square$

\begin{corollary}
	The sequence
	 $$\prod_{k=1}^N \frac 1{m_k} \Theta_k(z)\ov{\Theta_k(z)}\, d \dot z$$
	 of probabilistic  measures on the circle weakly converges.
\end{corollary}

This statement is well-known. To be complete, we give a proof.
Expressions $P_k(z)\,d\dot z$ are probabilistic measures,
since
the Fourier coefficient at $z^0$ is 1.  Other Fourier coefficients are increasing and 
$\le 1$. Therefore, we have a weak
convergence of measures.

\sm

{\bf\punct Proof of Theorem \ref{th:1}.}

\begin{lemma}
	The map $\cR$ is an isometry.
\end{lemma}

{\sc Proof.} It is sufficient to show that for two indicator
functions  $I[\sigmA]$, $I[\taU]$ we have
$$
\bigl\la \cR\, I_s[\sigmA], \cR\, I_t[\taU]\bigr\ra_{L^2(S^1,d\kappa)}=
\bigl\la I_s[\sigmA],  I_t[\taU] \bigr\ra_{L^2(\V)}=
\mu\bigl(B[\sigmA]\cap B[\taU]\bigr).
$$
Dividing balls into smaller subbals we reduce the statement to the case of balls
of the same level, say, $k$. In other words, we must prove the following identity
$$
\bigl\la
\cR\, U^p \,I_k[{0}_{\un k} \dots {0}_{\un 1};{0}],\,\,
\cR \, U^q\, I_k[{0}_{\un k} \dots {0}_{\un 1};{0}]\bigr\ra_{L^2(S^1,d\kappa)}=
\begin{cases}
\prod_{j=1}^k m_j^{-1},\quad \text{if $p=q$};
\\
0,\quad \text{otherwise},
\end{cases}
$$
where $p$, $q<h_k$.
The left hand side is
\begin{multline}
\lim_{N\to \infty}\int_{S^1} \frac{z^p}{\prod_{j=1}^k \Theta_j(z)}
\frac{\ov z^q}{\prod_{j=1}^k \ov{\Theta_j( z)}}\,\cdot \prod_{j=1}^N \frac 1{m_j} \Theta_j(z)\ov{\Theta_j(z)}\,d\dot z
=\\=
\prod_{j=1}^k \frac 1{m_j} \cdot
\lim_{N\to \infty}\int_{S^1} z^{p-q} \cdot \prod_{j=k+1}^N \frac 1{m_j} \Theta_j(z)\ov{\Theta_j(z)}\,d\dot z
\end{multline}
By Lemma \ref{eq:Riesz}.b, the prelimit integral is 1 if $p=q$.
If $p\ne q$, then by Lemma \ref{eq:Riesz}.c the integral is 0.
\hfill $\square$

\sm

Thus, the map $\cR$ extends to an isometric embedding of $L^2(\V)$ to
$L^2(S^1, d\kappa)$. 

\begin{lemma}
\label{l:last}
For any function $f\in L^2(\V)$, we have
$$
\cR\, U f(z)=z\cdot \cR f(z).  
$$
and 
$$
\cR\, U^{-1} f(z)=z^{-1}\cdot \cR f(z).  
$$
\end{lemma}

{\sc Proof.} Let us prove the first equality.
 It is sufficient to vary the statement
for functions $f=I_k[\sigmA]$. For 
$$
\sigmA\ne 
 \9_{\un k} \dots \9_{\un 2} \8_{\un 1}; \9_{\un 0}^{\un{1,m_1-2}}
$$
this is so by Corollary \ref{cor:}. For exceptional cases
we represent a ball $B_k[\dots]$ as union of subballs,
 The image of the last ball under the map $Q$ is a disjoint
union of balls of level $(k+1)$,
\begin{multline*}
 U\, I_k\bigl[\9_{\un k} \dots \9_{\un 2} \8_{\un 1};
 \9_{\un 0}^{{\un 1, \un {m_1-2}}}\bigr]=\\=
 U\,\Bigl(
I_{k+1}\bigl[0_{\un{k+1}}\,\9_{\un k} \dots \9_{\un 2} \8_{\un 1}; \9_{\un 0}^{{\un 1, \un {m_1-2}}}\bigr]
+ I_{k+1}\bigl[1_{\un{k+1}}\,\9_{\un k} \dots \9_{\un 2} \8_{\un 1}; \9_{\un 0}^{{\un 1, \un {m_1-2}}}\bigr]+
\dots \\
+
I_{k+1}\bigl[\9_{\un{k+1}}\9_{\un k} \dots
\9_{\un 2} \8_{\un 1}; \9_{\un 0}^{{\un 1, \un {m_1-2}}}\bigr]
\Bigr)
=\\= \sum_{i=0}^{\8_{k+1}}
I_{k+1}\bigl[i \9_{\un k} \dots \9_{\un 2} \9_{\un 1}; 0\bigr] +
U\, I_{k+1}\bigl[\9_{\un{k+1}} \9_{\un k} \dots \9_{\un 2} \8_{\un 1}; \9_{\un 0}^{{\un 1, \un {m_1-2}}}\bigr],
\end{multline*}
and we apply the same transformation to the last summand, etc.
In this way, we decompose the initial $I_k[\cdot]$ and 
$\cR I_k[\cdot]$ into two series consisting of disjoint indicator
functions. After this, we apply Corollary \ref{cor:}.a to each 
summand and pass to a limit.

The proof of the second statement is similar.
\hfill $\square$

\sm

{\sc End of proof of Theorem \ref{th:1}.}
The function $I_0[\mathbf{0}]$ is the indicator function 
of $\O\subset \V$. Its image under $\cR$ is $f(z)=1$.
By Lemma \ref{l:last}, $\cR$ sends $U^n I_0[\mathbf{0}]$ 
to $z^n$ for all $n\in\Z$.
 Therefore, the image 
of $\cR$ contains all polynomials in $z$, $z^{-1}$. By the Weierstrass
theorem, the image of $\cR$ is dense in the space of continuous functions on $S^1$
and therefore it is dense in $L^2(S^1,d\kappa)$.
\hfill $\square$


\section{An orthonormal  basis in $L^2(S^1,d\kappa)$}

\COUNTERS

{\bf\punct The construction of a basis.}
Consider a  system $\Phi[k,j,r;l, N](z)$ 
of rational functions of the circle defined in the following
way. The parameters
$k$, $j$, $r$, $l$, $N$ are nonnegative integers 
satisfying the conditions:

--- $k\ge 0$;

--- $0\le j\le m_k-2$;

--- $0\le r\le a_k^{j}$;

--- $l\ge k+1$;

--- if $l=k+1$, then $N=0$; otherwise $0<N< m_{k+2}\dots m_l$ 
and $m_l$ is not a divisor of $N$ .

The functions are defined by
\begin{multline*}
\Phi[k,j,r;l, N](z)=\\=
\frac{\sum\limits_{\taU=\tau_l\dots\tau_{k+2}: \tau_j<m_j}
\exp\Bigl\{\frac{2\pi i}{m_{k+2}\dots m_l}\cdot N\cdot \ov{\ov{\tau_l\dots\tau_{k+2}}} \Bigr\}
z^{\Upsilon_l( \tau_l\dots\tau_{k+2}j\,\9_k\dots \9_1;r)}}
{\prod_{s=1}^l \Theta_s(z)},
\end{multline*}
where $\ov{\ov \taU}$ is defined by \eqref{eq:ov} and $\Upsilon_l(\cdot)$
by \eqref{eq:Upsilon}.

\begin{proposition}
\label{pr:}
The system $\Phi[k,j,r; l, N](z)$
is an orthogonal basis in $L^2(S^1,d\kappa)$,
$$
\|\Phi[k,j,r; l,N]\|^2_{L^2}=\frac 1{m_{k+2}m_{k+3}\dots m_{1}}.
$$ 
\end{proposition}

\sm

{\bf\punct Proof of Proposition \ref{pr:}.}
Consider an additive group of the ring $\O[n_1,n_2,\dots]$.
It is a compact Abelian group and it is the inverse limit
$$
\lim\limits_{\longleftarrow} \Z/n_1\dots n_k\Z.
$$
Therefore the Pontryagin dual group
$\O(n_1,n_2,\dots)^\wedge$ is the direct limit
$$
\lim\limits_{\longrightarrow} (\Z/n_1\dots n_k\Z)^\wedge.
$$
This group is discrete, its elements (characters of $\O$) form
an orthonormal basis in $L^2\bigl(\O[n_1,n_2,\dots] \bigr)$ 
Since $(\Z/n_1\dots n_k\Z)^\wedge$ is isomorphic to the group   
$\Z/n_1\dots n_k\Z$ itself.
The characters of $\Z/n_1\dots n_k\Z$ are 
given by the formula
$$
\chi_\alpha(l)=
\exp\Bigl\{\frac{2\pi i \,\alpha\, l}{n_1\dots n_k}  \Bigr\},
$$
both $\alpha$ and $l$  are defined modulo $n_1\dots n_k$.
If $n_k$ is a divisor of $\alpha$, then actually 
$\chi_\alpha$ is defined on the quotient group $\Z/n_1\dots n_{k-1}\Z$.
So the set of characters is a disjoint union
$$
\coprod_k \Bigl((\Z/n_1\dots n_k\Z)^\wedge\setminus (\Z/n_1\dots n_{k-1}\Z)^\wedge\Bigr)
$$
Thus, each character of $\O(n_1,n_2,\dots)$ has
the form 
$$
\sum_{\sigmA=\sigma_k\dots \sigma_1:\, \sigma_j<n_j}
\exp\Bigl\{\frac{2\pi i}{n_1\dots n_k}\, N\,\ov{\ov{\sigma_k\dots \sigma_1}}  \Bigr\}
I_k[\sigma_k\dots \sigma_1]
$$
for some $k$ and $N$ such that $n_k$ is not a divisor of $N$

Now we represent the space $\V$ as a disjoint union 
of balls
$$
\V=\coprod_{k,j,r}B_{k+1}[j\,\9_k\dots\9_1;r],
$$
where $j\ne m_{k+1}-1$. Such a ball
has measure $1/m_1\dots m_{k+1}$ and  is in the obvious one-to-one 
correspondence with the set $\O(m_{k+2},m_{k+3},\dots)$.
We take the system of characters on each ball and get an orthogonal
basis in $L^2(\V)$.

It remains to apply our transformation $L^2(\V)\to L^2(S^1,d\kappa)$ to characters of  
$\O(m_{k+2},m_{k+3},\dots)$.

\sm

{\bf Acknowledgements.} I am grateful to V.~V.~Ryzhikov, M.~E.~Lipatov,
and M.~S.~Lobanov for discussions of this topics.

 \tt
 
University of Graz,\\
\vphantom{.}\hfill Department of Mathematics and Scientific Computing;

High School of Modern Mathematics MIPT; 

Moscow State University, MechMath Faculty;

Univesity of Vienna, Faculty of Mathematics.
 
 \noindent

 URL: http://mat.univie.ac.at/$\sim$neretin/

\end{document}